\documentclass[24pt]{amsart}
\usepackage{amsmath, amsthm, amscd, amsfonts}

\newtheorem{thm}{Theorem}[section]

\newtheorem{cor}[thm]{Corollary}
\theoremstyle{definition}

\theoremstyle{remark}

\numberwithin{equation}{section}
\newfont{\kh}{msbm10}

\def\N{{\mathbb N}}

\def\S{\mathfrak{S}}

\begin{document}

\title[On the decomposition of $n!$ into primes]{On the decomposition of $n!$ into primes}

\author{Mehdi Hassani}

\address{Mehdi Hassani, \newline Department of Mathematics, Institute for Advanced Studies in Basic
Sciences, P.O. Box 45195-1159, Zanjan, Iran}

\email{mmhassany@member.ams.org}



\begin{abstract}
In this note, we make explicit approximation of the average of prime
powers in the decomposition of $n!$. Also, some related questions
have been studied.
\end{abstract}
\maketitle

\section{Introduction}
Letting
$$
n!=\prod_{p\leq n}p^{v_p(n!)},
$$
with $p$ is prime, it is known \cite{nathan}, as a classic result
that
\begin{equation}\label{vpsum}
v_p(n!)=\sum_{k=1}^\infty\left\lfloor\frac{n}{p^k}\right\rfloor=\sum_{k=1}^m\left\lfloor\frac{n}{p^k}\right\rfloor,
\end{equation}
with $m=m_{n,p}=\lfloor\frac{\log n}{\log p}\rfloor$ and $\lfloor
x\rfloor$ is the largest integer less than or equal to $x$. In
this paper, we study the following summation for a fixed positive
integer $n$,
$$
\Upsilon(n)=\sum_{p\leq n}v_p(n!).
$$
\subsection{Approximate Formula for the Function $\Upsilon(n)$}
First, we note that integrating by parts, yields
\begin{eqnarray}\label{int1log}
\int_2^n\frac{dx}{\log x}&=&n\sum_{k=1}^N\frac{(k-1)!}{\log^k
n}-2\sum_{k=1}^N\frac{(k-1)!}{\log^k
2}+N!\int_2^n\frac{dx}{\log^{N+1}x}\\\nonumber
&=&n\sum_{k=1}^N\frac{(k-1)!}{\log^k
n}+O\Big(\frac{n}{\log^{N+1}n}\Big).
\end{eqnarray}
Considering (\ref{vpsum}), we have
$$
\Upsilon(n)=\sum_{p\leq n}\sum_{k\leq
m}\left\lfloor\frac{n}{p^k}\right\rfloor=\sum_{p\leq n}\sum_{k\leq
m}\Big(\frac{n}{p^k}+O(1)\Big).
$$
So, we have
\begin{eqnarray*}
\Upsilon(n)-\sum_{p\leq n}\sum_{k\leq
m}\frac{n}{p^k}&\ll&\sum_{p\leq n}\sum_{k\leq m}1\ll\sum_{p\leq
n}m\ll\log n\sum_{p\leq n}\frac{1}{\log p}\\&<&\log n\sum_{k\leq
n}\frac{1}{\log k}\ll\log n\int_2^n\frac{dx}{\log x},
\end{eqnarray*}
and using (\ref{int1log}) with $N=1$, we obtain
$$
\Upsilon(n)-\sum_{p\leq n}\sum_{k\leq m}\frac{n}{p^k}\ll n.
$$
Thus, since $m\geq 1$, we have
$$
\Upsilon(n)=n\sum_{p\leq n}\sum_{k\leq
m}\frac{1}{p^k}+O(n)=n\sum_{p\leq
n}\frac{1-\frac{1}{p^m}}{p-1}+O(n)=n\sum_{p\leq
n}\frac{1}{p}+O(n).
$$
In the other hand, it is known \cite{apos} that
$$
\sum_{p\leq n}\frac{1}{p}=\log\log n+O(1).
$$
Therefore,
$$
\Upsilon(n)=n\log\log n+O(n).
$$
Now, let $\overline{\Upsilon}(n)$ be the mean value of the values
of $v_p(n!)$ for $p\leq n$. We have
$$
\overline{\Upsilon}(n)=\frac{1}{\#\{v_p(n!)|p\leq n\}}\sum_{p\leq
n}v_p(n!)=\frac{\Upsilon(n)}{\pi(n)},
$$
where $\pi(n)=$ the number of primes not exceeding of $n$.
Considering the Prime Number Theorem (PNT) \cite{apos};
$\pi(n)\sim\frac{n}{\log n}$, we obtain
$$
\overline{\Upsilon}(n)=\frac{n\log\log
n}{\pi(n)}+O\Big(\frac{n}{\pi(n)}\Big)=\log n\log\log n+O(\log n).
$$
What does this mean? Putting $\mathfrak{L}=\log n$ and letting
$p_x=\lfloor x\rfloor^{th}$ prime number for $x\geq 1$, another
analogue of PNT yields that $\overline{\Upsilon}(n)\sim
p_{\mathfrak{L}}$, which means the average of the prime powers in
the factorization of $n!$ into the primes is approximately
$\mathfrak{L}^{th}$ prime number.

\subsection{Aim of Work and Summary of the Results} In the next sections,
first we get some explicit bounds for the function $\Upsilon(n)$,
and then consequently for the function $\overline{\Upsilon}(n)$.
More precisely, we prove the following results. Note that the
constants $c_4$, $c_8$ and $c_{10}$ at bellow all are effective.
\begin{thm}\label{eub-g} For every $n\geq 2$, we have
$$
\Upsilon(n)<(n-1)\log\log(n-1)+c_4(n-1)+\frac{n}{\log
n}+\frac{1717433 n}{\log^5 n}.
$$
\end{thm}

\begin{thm}\label{eub-gbar} For every $n\geq 3$, we have
$$
\overline{\Upsilon}(n)<\frac{\log n}{1+\log n}\log n\log\log
(n-1)+\frac{c_4\log^2 n}{1+\log n}+\frac{\log n}{1+\log
n}+\frac{1717433}{(1+\log n)\log^3 n}.
$$
\end{thm}

\begin{cor}\label{eub-gbar-cor} For $n\geq 12602987$, we have
$$
\overline{\Upsilon}(n)<\log n \log\log n +{ \frac
{380537}{17966}}\log n +1.
$$
\end{cor}

\begin{thm}\label{elb-g} For every $n\geq 3$, we have
$$
\Upsilon(n)>(n-1)\log\log n+c_8(n-1)-\frac{n}{\log
n}-\frac{16381n}{5000\log^2 n}-\frac{6n}{\log^3
n}-\frac{54281n}{800\log^4 n}-c_{10}\log n.
$$
\end{thm}

\begin{thm}\label{elb-gbar} For every $n\geq 2$, we have
\begin{eqnarray*}
\overline{\Upsilon}(n)&>&\frac{(n-1)\kappa_n}{n}\log n\log\log
n+\frac{c_8(n-1)\kappa_n\log n}{n}-\frac{16381\kappa_n}{5000\log
n}-\frac{6\kappa_n}{\log^2 n}\\&-&\frac{54281\kappa_n}{800\log^3
n}-\frac{c_{10}\kappa_n\log^2 n}{n},
\end{eqnarray*}
where
$$
\kappa_n=\frac{5000\log n}{6381+5000\log n}.
$$
\end{thm}

\subsection{Some Tools}
During proofs, we will need to estimate summations of the form
$\sum_{p\leq n}f(p)$ for a given function $f(x)\in
C^1(\mathbb{R}^+)$ with summation over primes $p$. Concerning this
problem, using Stieljes integral \cite{r-s-62} and integrating by
parts, we have
\begin{eqnarray}\label{sumpleqn}
\sum_{p\leq n}f(p)=\int_{2^-}^n \frac{f(x)}{\log
x}d\vartheta(x)=\frac{f(n)\vartheta(n)}{\log
n}+\int_2^n\vartheta(x)\frac{d}{dx}\left(\frac{-f(x)}{\log
x}\right)dx,
\end{eqnarray}
where $\vartheta(x)=\sum\limits_{p\leq x}\log p$, and it is known
that \cite{dusart} for $x>1$, we have
\begin{equation}\label{tb2}
|\vartheta(x)-x|<\frac{793x}{200\log^2x},
\end{equation}
and
\begin{equation}\label{tb4}
|\vartheta(x)-x|<1717433\frac{x}{\log^4x}.
\end{equation}
Starting point of explicit approximations of $\Upsilon(n)$ is the
following known \cite{hassani} bounds
\begin{equation}\label{vpbns}
\frac{n-p}{p-1}-\frac{\log n}{\log p}<v_p(n!)\leq \frac{n-1}{p-1},
\end{equation}
which holds true for every $n\in\N$ and prime $p$, with $p\leq n$.
To apply obtained results for approximating
$\overline{\Upsilon}(n)$, we need some explicit bounds concerning
$\pi(n)$; it is known \cite{dusart} that
\begin{equation}\label{pi-lb}
\pi(n)\geq \frac{n}{\log n}\left(1+\frac{1}{\log n}\right),
\end{equation}
which holds true for every $n\geq 599$. Also, for every $n\geq 2$,
we have
\begin{equation}\label{pi-ub}
\pi(n)\leq \frac{n}{\log n}\left(1+\frac{6381}{5000\log n}\right).
\end{equation}
To do careful computations, we use the Maple software. Specially,
to compute the values of $\Upsilon(n)$ (and consequently
$\overline{\Upsilon}(n)$), we use the following program in Maple
software worksheet:
\\
\\
\textsf{ G:=proc(n)\\
tot := 0:\\
for i from 1 by 1 while ithprime(i)$<$n do\\
tot := tot +
sum(floor(n/ithprime(i)**k),k=1..floor(log(n)/log(ithprime(i))))\\
end do:\\
end:}
\\

\section{Explicit Approximation of the Functions $\Upsilon(n)$ and $\overline{\Upsilon}(n)$}
In this section we introduce the proof of mentioned explicit
bounds for the functions $\Upsilon(n)$ and
$\overline{\Upsilon}(n)$.
\subsection{Upper Bounds} Using the right hand
side of (\ref{vpbns}) and (\ref{sumpleqn}), we have
$\Upsilon(n)\leq S_1(n)$, where
\begin{eqnarray}\label{s1}
S_1(n)&=&\sum_{p\leq n}\frac{n-1}{p-1}\\\nonumber
&=&\frac{\vartheta(n)}{\log
n}+(n-1)\int_{2}^n\vartheta(x)\frac{d}{dx}\left(\frac{-1}{(x-1)\log
x}\right)dx.
\end{eqnarray}
\subsubsection{Upper Approximation of $S_1(n)$} Since,
$\frac{d}{dx}(\frac{-1}{(x-1)\log x})>0$, using (\ref{tb2}), we
obtain
$$
\int_{2}^n\vartheta(x)\frac{d}{dx}\left(\frac{-1}{(x-1)\log
x}\right)dx<\mathcal{I}_1(n)+\mathcal{E}_1(n)+c_1,
$$
where
$$
\mathcal{I}_1(n)=\int_{2}^n \\{\frac
{1200\,{x}^{3}+365\,{x}^{2}+9944\,x- 1993}{1200 \left( x-1 \right)
^{4}\log x}}{dx},
$$
and
$$
c_1=\frac{-5937\log^2 2+3965\log 2+1586}{600\log^3 2}\approx
7.416262921,
$$
and $\mathcal{E}_1(n)=-\frac {A(n)}{B(n)}$ with
$B(n)=1200(n-1)^3\log^3 n$, and
\begin{eqnarray*}
\frac{A(n)}{n}&=&1200 n^2\log^2 n+2379\,{n}^{2}\log n+1586
n^2-6365 n \log^2 n\\&-&3172 n\log n-3172n+1993\log^2 n+793\log
n+1586.
\end{eqnarray*}
Easily $\lim\limits_{n\rightarrow\infty}\mathcal{E}_1(n)\log n=-1$
and for every $n$ we have $\mathcal{E}_1(n)<0$. Therefore, we get
$$
S_1(n)<\frac{\vartheta(n)}{\log
n}+c_1(n-1)+(n-1)\mathcal{I}_1(n)\hspace{10mm}(n\geq 2).
$$
Now, we have
\begin{eqnarray*}
\mathcal{I}_1(n)=\int_{e+1}^n {\frac
{1200\,{x}^{3}+365\,{x}^{2}+9944\,x- 1993}{1200 \left( x-1 \right)
^{4}\log x}}{dx}+c_2,
\end{eqnarray*}
where
$$
c_2=\int_2^{e+1} {\frac {1200\,{x}^{3}+365\,{x}^{2}+9944\,x-
1993}{ 1200\left( x-1 \right) ^{4}\log x }}{dx}\approx
12.35466367,
$$
and so,
$$
\mathcal{I}_1(n)<\int_{e+1}^n {\frac
{1200\,{x}^{3}+365\,{x}^{2}+9944\,x- 1993}{1200 \left( x-1 \right)
^{4}\log(x-1) }}{dx}+c_2=\log\log(n-1)+\mathcal{E}_2(n)+c_3,
$$
where
\begin{eqnarray*}
\mathcal{E}_2(n)=&-&{\frac {793}{240}}\,{\it Ei} \left( 1,\log
\left( n-1 \right)
 \right) -{\frac {2379}{200}}\,{\it Ei} \left( 1,2\,\log  \left( n-1
 \right)  \right)\\ &-&{\frac {793}{100}}\,{\it Ei} \left( 1,3\,\log
 \left( n-1 \right)  \right)\rightarrow 0^-\hspace{10mm}(n\geq 2),
\end{eqnarray*}
and
$$
c_3={\frac {793}{240}}\,{\it Ei} \left( 1,1 \right) +{\frac
{2379}{200}}\, {\it Ei} \left( 1,2 \right) +{\frac
{793}{100}}\,{\it Ei} \left( 1,3
 \right) + c2\approx 13.76468999.
$$
Note that $Ei$ is the formal notation for the Exponential Integral
\cite{a-s}, defined by
$$
{\it Ei} \left( a,z \right) =\int _{1}^{\infty
}\!{e^{-tz}}{t}^{-a}{dt}\hspace{10mm}(\Re(z)>0).
$$
Therefore, putting $c_4=c_1+c_3\approx 21.18095291$, we obtain
$$
S_1(n)<(n-1)\log\log(n-1)+c_4(n-1)+\frac{\vartheta(n)}{\log
n}\hspace{10mm}(n\geq 2),
$$
and using (\ref{tb4}), we get the following explicit upper bound
$$
S_1(n)<(n-1)\log\log(n-1)+c_4(n-1)+\frac{n}{\log n}+\frac{1717433
n}{\log^5 n}\hspace{10mm}(n\geq 2).
$$
Remembering $\Upsilon(n)\leq S_1(n)$, completes the proof of the
Theorem \ref{eub-g}. Now, we can use this result to get some upper
bounds for the function $\overline{\Upsilon}(n)$. Since
$\overline{\Upsilon}(n)=\frac{\Upsilon(n)}{\pi(n)}$, considering
(\ref{pi-lb}), for every $n\geq 599$ we have
$$
\overline{\Upsilon}(n)<\frac{\log n}{1+\log n}\log n\log\log
(n-1)+\frac{c_4\log^2 n}{1+\log n}+\frac{\log n}{1+\log
n}+\frac{1717433}{(1+\log n)\log^3 n},
$$
which holds true for $3\leq n\leq 598$ too, by computation. This
proofs the Theorem \ref{eub-gbar}. Also, an straight computation
yields the following simpler bound for $n\geq 12602987$,
$$
\overline{\Upsilon}(n)<\log n \log\log n +{ \frac
{380537}{17966}}\log n +1.
$$
This proofs the Corollary \ref{eub-gbar-cor}.

\subsection{Lower Bounds}

Using the right hand side of (\ref{vpbns}) and (\ref{sumpleqn}),
we have
\begin{eqnarray}\label{uplb}
\Upsilon(n)>\sum_{p\leq n}\left(\frac{n-p}{p-1}-\frac{\log n}{\log
p}\right)=S_1(n)-\pi(n)-S_2(n),
\end{eqnarray}
where $S_1(n)$ has been introduced in (\ref{s1}), and
\begin{eqnarray}\label{s2}
S_2(n)=\sum_{p\leq n}\frac{\log n}{\log
p}=\frac{\vartheta(n)}{\log n}+\log
n\int_{2}^n\vartheta(x)\frac{d}{dx}\left(\frac{-1}{\log^2
x}\right)dx.
\end{eqnarray}
\subsubsection{Lower Approximation of $S_1(n)$.} Because
$\frac{d}{dx}\left(\frac{-1}{(x-1)\log x}\right)>0$, considering
(\ref{tb2}), we have
$$
\int_{2}^n\vartheta(x)\frac{d}{dx}\left(\frac{-1}{(x-1)\log
x}\right)dx>\mathcal{I}_2(n)+\mathcal{E}_3(n)+c_5,
$$
where
$$
\mathcal{I}_2(n)=\int_{2}^n \\{\frac
{1200x^3-7565x^2-2744x-407}{1200 (x-1)^{4}\log x}}{dx},
$$
and
$$
c_5=\frac{8337\log^2 2-3965\log 2-1586}{600\log^3 2}\approx
-1.645482755,
$$
and $\mathcal{E}_3(n)=-\frac {C(n)}{D(n)}$ with
$D(n)=1200(n-1)^3\log^3 n$, and
\begin{eqnarray*}
\frac{C(n)}{n}&=&1200 n^2\log^2 n-2379\,{n}^{2}\log n-1586
n^2+1565 n \log^2 n\\&+&3172 n\log n+3172n+407\log^2 n-793\log
n-1586.
\end{eqnarray*}
Easily $\lim\limits_{n\rightarrow\infty}\mathcal{E}_3(n)\log
n=-1$. The function $\mathcal{E}_3(n)$ takes its minimum vale at
$n\approx 28.85589912$. Thus for every $n\geq 2$, we have
\begin{eqnarray*}
\mathcal{E}_3(n)&>&\min\{\mathcal{E}_3(28),\mathcal{E}_3(29)\}=\mathcal{E}_3(29)\\
&=&-\frac{29(131874\log^2 29-238693\log 29-155428)}{3292800\log^3
29}\approx -.1236613745.
\end{eqnarray*}
In the other hand, we have
$$
\mathcal{I}_2(n)=\int_{e}^n \\{\frac
{1200x^3-7565x^2-2744x-407}{1200 (x-1)^{4}\log x}}{dx}+c_6,
$$
where
$$
c_6=\int_{2}^e \\{\frac {1200x^3-7565x^2-2744x-407}{1200
(x-1)^{4}\log x}}{dx}\approx -8.600279758.
$$
So,
$$
\mathcal{I}_2(n)>\int_{e}^n \\{\frac
{1200x^3-7565x^2-2744x-407}{1200 x^{4}\log x}}{dx}+c_6=\log\log
n+\mathcal{E}_4(n)+c_7,
$$
where
$$
\mathcal{E}_4(n)={\frac {1513}{240}}\,{\it Ei} \left( 1,\log n
\right) +{\frac {343}{150}}\,{\it Ei} \left( 1,2\,\log n
 \right) +{\frac {407}{1200}}\,{\it Ei} \left( 1,3\,\log n
\right),
$$
and
$$
c_7=c_6-\left({\frac {1513}{240}}\,{\it Ei} \left( 1,1\right)
+{\frac {343}{150}}\,{\it Ei} \left( 1,2 \right) +{\frac
{407}{1200}}\,{\it Ei} \left( 1,3\right)\right)\approx
-10.09955739.
$$
Note that, $\frac{d}{dn}\mathcal{E}_4(n)=-\big({\frac
{1513}{240{n}^{2}\log n }}+{\frac {343}{150{n}^{3}\log n}}+{\frac
{407}{1200{n}^{4}\log n}}\big)<0$ and
$\lim\limits_{n\rightarrow\infty}\mathcal{E}_4(n)=0$. Thus, for
every $n\geq 2$, we have $\mathcal{E}_4(n)>0$. Therefore, we
obtain
\begin{equation}\label{s1lbt}
S_1(n)>\frac{\vartheta(n)}{\log n}+(n-1)\log\log n+c_8(n-1),
\end{equation}
where $c_8=c_5+c_7+\mathcal{E}_3(29)\approx -11.86870152$.
Considering (\ref{tb4}), we get the following explicit lower bound
for every $n\geq 2$
$$
S_1(n)>(n-1)\log\log n+c_8(n-1)+\frac{n}{\log
n}-\frac{1717433n}{\log^5 n}.
$$
\subsubsection{Lower Approximation of $S_2(n)$.} Because
$\frac{d}{dx}\left(\frac{-1}{\log^2 x}\right)>0$, considering
(\ref{tb2}), we have
\begin{eqnarray*}
\int_{2}^n\vartheta(x)\frac{d}{dx}\left(\frac{-1}{\log^2
x}\right)dx>\int_{2}^n\frac{200\log^2 x-793}{100\log^5
x}dx=\frac{1607}{2400}\int_2^n\frac{dx}{\log
x}+\mathcal{R}_1(n)+c_9,
\end{eqnarray*}
where
$$
\mathcal{R}_1(n)=-\frac {1607n}{2400\log
n}-\frac{1607n}{2400\log^2 n}+\frac {793n}{1200\log^3
n}+\frac{793n}{400\log^4 n},
$$
and $c_9=-\mathcal{R}_1(2)\approx -16.42613005$. Now, considering
(\ref{int1log}), and a simple calculation, yields that
$$
\int_2^n\frac{dx}{\log x}>n\sum_{k=1}^{5}\frac{(k-1)!}{\log^k
x}\hspace{10mm}(n\geq 563.74).
$$
Applying this bound, we obtain
\begin{eqnarray*}
\int_{2}^n\vartheta(x)\frac{d}{dx}\left(\frac{-1}{\log^2
x}\right)dx>\frac{2n}{\log^3 n}+\frac{6n}{\log^4
n}+\frac{1607n}{100\log^5 n}+c_9.
\end{eqnarray*}
Therefore,
\begin{equation*}
S_2(n)>\frac{\vartheta(n)}{\log n}+\frac{2n}{\log^2
n}+\frac{6n}{\log^3 n}+\frac{1607n}{100\log^4 n}+c_9\log
n\hspace{10mm}(n\geq 564),
\end{equation*}
and considering (\ref{tb4}), we obtain
$$
S_2(n)>\frac{n}{\log n}+\frac{2n}{\log^2 n}+\frac{6n}{\log^3
n}+\frac{1607n}{100\log^4 n}-\frac{1717433n}{\log^5 n}+c_9\log
n\hspace{10mm}(n\geq 564).
$$
\subsubsection{Upper Approximation of $S_2(n)$.} Again, considering the
relations $\frac{d}{dx}\left(\frac{-1}{\log^2 x}\right)>0$ and
(\ref{tb2}), we have
\begin{eqnarray*}
\int_{2}^n\vartheta(x)\frac{d}{dx}\left(\frac{-1}{\log^2
x}\right)dx<\int_{2}^n\frac{200\log^2 x+793}{100\log^5
x}dx=\frac{3193}{2400}\int_2^n\frac{dx}{\log
x}+\mathcal{R}_2(n)+c_{10},
\end{eqnarray*}
where
$$
\mathcal{R}_2(n)=-\frac {3193n}{2400\log
n}-\frac{3193n}{2400\log^2 n}-\frac {793n}{1200\log^3
n}-\frac{793n}{400\log^4 n},
$$
and $c_{10}=-\mathcal{R}_2(2)\approx 30.52238614$. Now, an easy
computation yields that
$$
\varepsilon+\int_2^n\frac{dx}{\log
x}<n\sum_{k=1}^{4}\frac{(k-1)!}{\log^k
x}+\frac{51n}{\log^5n}\hspace{10mm}(n\geq 2~{\rm and}~
\varepsilon\approx 0.144266447).
$$
Thus, for every $n\geq 2$ we have
\begin{eqnarray*}
\int_{2}^n\vartheta(x)\frac{d}{dx}\left(\frac{-1}{\log^2
x}\right)dx<\frac{2n}{\log^3 n}+\frac{6n}{\log^4
n}+\frac{54281n}{800\log^5 n}+c_{10}.
\end{eqnarray*}
Therefore,
\begin{equation}\label{s2ubt}
S_2(n)<\frac{\vartheta(n)}{\log n}+\frac{2n}{\log^2
n}+\frac{6n}{\log^3 n}+\frac{54281n}{800\log^4 n}+c_{10}\log n,
\end{equation}
and considering (\ref{tb4}), we obtain
$$
S_2(n)<\frac{n}{\log n}+\frac{2n}{\log^2 n}+\frac{6n}{\log^3
n}+\frac{54281n}{800\log^4 n}+\frac{1717433n}{\log^5 n}+c_{10}\log
n.
$$
Therefore, considering the relations (\ref{uplb}), (\ref{s1lbt})
and (\ref{s2ubt}), for every $n\geq 2$ we obtain
$$
\Upsilon(n)>-\pi(n)-\frac{2n}{\log^2 n}-\frac{6n}{\log^3
n}-\frac{54281n}{800\log^4 n}+(n-1)\log\log n+c_8(n-1)-c_{10}\log
n,
$$
and considering (\ref{pi-ub}), we get
$$
\Upsilon(n)>(n-1)\log\log n+c_8(n-1)-\frac{n}{\log
n}-\frac{16381n}{5000\log^2 n}-\frac{6n}{\log^3
n}-\frac{54281n}{800\log^4 n}-c_{10}\log n.
$$
This completes the proof of the Theorem \ref{elb-g}. Dividing both
sides of above inequality by $\pi(n)$ and using (\ref{pi-ub}), we
obtain
\begin{eqnarray*}
\overline{\Upsilon}(n)&>&\frac{(n-1)\kappa_n}{n}\log n\log\log
n+\frac{c_8(n-1)\kappa_n\log n}{n}-\frac{16381\kappa_n}{5000\log
n}-\frac{6\kappa_n}{\log^2 n}\\&-&\frac{54281\kappa_n}{800\log^3
n}-\frac{c_{10}\kappa_n\log^2 n}{n},
\end{eqnarray*}
where
$$
\kappa_n=\frac{5000\log n}{6381+5000\log n}.
$$
This gives the proof of the Theorem \ref{elb-gbar}.

\section{Some Questions and Answers}

\subsection{Approximately, at which prime $\overline{\Upsilon}(n)$ appear?} To answer this, we have to solve the equation
$\overline{\Upsilon}(n)=v_p(n!)$ approximately, according to $p$.
Using the relation (\ref{vpbns}), we have
\begin{equation}\label{vpno}
v_p(n!)=\frac{n}{p-1}+O(\log n).
\end{equation}
Putting this and the relation $\overline{\Upsilon}(n)=\log
n\log\log n+O(\log n)$ in the approximate equation
$\overline{\Upsilon}(n)=v_p(n!)$, we obtain
$$
p=\frac{n}{\log n\log\log n}+1+O\Big(\frac{1}{\log
n}\Big)\sim\frac{n}{\log n\log\log
n}\hspace{10mm}(n\rightarrow\infty).
$$
If we let $p$ to be (approximately) the $k^{th}$ prime, then
considering PNT we have
$$
\frac{n}{\log n\log\log n}\sim k\log
k\hspace{10mm}(n\rightarrow\infty).
$$
Solving this approximate equation according to $k$, we obtain
$$
k\sim\frac{n}{\log n\log\log n~W\big(\frac{n}{\log n\log\log
n}\big)}\hspace{10mm}(n\rightarrow\infty).
$$
where $W$ is the Lambert $W$ function, defined by $W(x)e^{W(x)}=x$
for $x\in[-e^{-1},+\infty)$, and it known \cite{cghjk} that
$W(x)\sim\log x$ when $x\rightarrow\infty$. Therefore
$$
k\sim\frac{n}{\log n\log\log n\log\big(\frac{n}{\log n\log\log
n}\big)}\sim\frac{n}{\log^2 n\log\log
n}\hspace{10mm}(n\rightarrow\infty).
$$
This means that $\overline{\Upsilon}(n)$ appears approximately at
$k^{th}$ prime with above obtained $k$.

\subsection{Minimal Square Perfecter of $n!$} We define (as usual)
the minimal perfecter of the positive integer $N$ as follows
$$
\S(n)=\min\Big\{m|m\in\N,~mn=l^2~{\rm for~some}~l\in\N\Big\}.
$$
It is clear that $1\leq\mathfrak{S}(n!)\leq \prod\limits_{p\leq
n}p=e^{\vartheta(n)}$. To clarify situation, we define
$$
\vartheta(n;q,a)=\sum_{\substack{
p\le n\\
v_p(n!)\equiv a\pmod q}} \log p,
$$
for a given prime $q$, if there exists some $p$ with $v_p(n!)\equiv
a\pmod q$, and else by putting 1. Note that with this definition we
have $\vartheta(n)=\sum\limits_{0\leq a\leq q-1}\vartheta(n;q,a)$.
An easy computation implies that for $n>1$, the Bertrand's
conjecture is equivalent by $1<\mathfrak{S}(n!)$. However, we can
find more better lower bounds; for $n\geq 4$ we have
$$
\S(n!)=e^{\vartheta(n;2,1)}\geq e^{\left(\sum\limits_{v_p(n!)=1}\log
p\right)}=e^{\left(\sum\limits_{\sqrt{n}<p\leq n<2p}\log
p\right)}=e^{\left(\sum\limits_{\frac{n}{2}<p\leq n}\log
p\right)}=e^{\vartheta(n)-\vartheta(\frac{n}{2})}.
$$
Therefore, considering (\ref{tb2}), we get the following bounds for
$n\geq 4$,
$$
e^{\frac{n}{2}-\frac{793}{200}n\big(\frac{1}{\log
n}+\frac{1}{2\log(\frac{n}{2})}\big)}<\S(n!)<e^{n+\frac{793n}{200\log
n}}.
$$
Here two questions arise to mind:\\
\textit{Question 1.} How we can estimate $\vartheta(n;q,a)$? To
answer this, we need an equivalent condition by $v_p(n!)\equiv
a\pmod q$.\\
\textit{Question 2.} Insomuch, Bertrand's conjecture is equivalent
by $1<\mathfrak{S}(n!)$, what we can get about $\S(n!)$ under
assumption recent improvements concerning existing primes in short
intervals?\\
A very good result on the existing primes in short intervals is due
to Baker, Harman and Pintz \cite{b-h-p} which asserts that there
exists real $x_0$ such that for all $x>x_0$ the interval
$\big[x-x^{0.525},x\big]$ contains a prime.\\
\textit{Question 3.} What is exact value of $\mathfrak{S}(n!)$?

\end{document}